\documentclass{commat}


\newcommand{\F}{\mathcal{F}}

\title{%
    On the geometric mean of the values of positive multiplicative arithmetical functions
    }

\author{%
    Mehdi Hassani and Mohammadreza Esfandiari
    }

\affiliation{
    \address{
    Department of Mathematics
    University of Zanjan
    University Blvd., 45371-38791, Zanjan
    Iran
    }
    \email{%
    mehdi.hassani@znu.ac.ir, 
    m.esfandiari@znu.ac.ir
    }
}

\abstract{%
   In this paper we obtain asymptotic expansions for the geometric mean of the values of positive strongly multiplicative function $f$ satisfying \[
   f(p) = \alpha (d) p^d + O(p^{d-\delta})
   \]
   for any prime $p$, with $d$ real, $\alpha (d)$ and $\delta>0$.
    }



\VOLUME{31}
\YEAR{2023}
\NUMBER{1}
\firstpage{21}
\DOI{https://doi.org/10.46298/cm.10133}

\begin{paper}

\section{Introduction}

The arithmetic mean of the values of arithmetical functions is very well studied in the literature and some theories have been developed. For example, see \cite{dek-luca}, \cite{elliott-85}, \cite{elliott-97}, \cite{elliott-79}, \cite{elliott-80}, \cite{iwa-kowal}, \cite{kuku}, \cite{tenen} and the references given there. In comparison, the geometric mean of the values of arithmetical functions has been studied only in some special cases. Let us denote by $G_f(n)$ the geometric mean of the first $n$ values of the positive arithmetic function $f$. In 2008 Deshouillers and Luca \cite{DL} studied the density modulo 1 of some sequences involving the values of the Euler function $\varphi$, including the sequence with general term $G_\varphi(n)$. Meanwhile, they proved that
\[
G_\varphi(n)=\frac{1}{\mathrm{e}}\prod_{p}\left(1-\frac{1}{p}\right)^\frac{1}{p}n+O(\log n).
\]
In 2013 the first author \cite{MH-ud} studied uniform distribution modulo 1 of some sequences involving the values of the Euler function, where he improved the above error term up to $O(\log\log n)$. In 2012 Bosma and Kane \cite{bosma-kane}, and later in 2018 Pomerance \cite{pomerance} considered a variant of the geometric mean of the first $n$ values of the sum of divisors function to study the so-called ``aliquot constant''. In 2016 the first author \cite{hass-bims} considered the geometric mean of the first $n$ values of the function $\mathrm{d}(n)=\sum_{d|n}1$. He proved that given any positive integer $r$, there exist computable constants $c_1,\dots,c_r$ such that
\[
G_\mathrm{d}(n)=2^M\prod_{\substack{p^\alpha\\ \alpha\geqslant 2}}\log\Big(1+\frac{1}{\alpha}\Big)^\frac{1}{p^\alpha}\left(\log n\right)^{\log 2}\left(1+\sum_{j=1}^r \frac{c_{j}}{\log^{j}n}+O\Big(\frac{1}{\log^{r+1}n}\Big)\right),
\]
where $M$ is the Meissel--Mertens constant. Recently, the second author \cite{esfandiari} studied the geometric mean of the first $n$ values of the Jordan totient function.

In this paper we are motivated by introducing a general theory for studying the geometric mean of the values of positive multiplicative functions $f$. Since $\log f$ is additive, we get
\[
n\log G_f(n)
=\sum_{k\leqslant n}\log f(k)
=\sum_{\substack{p^\alpha\leqslant n\\ \alpha\geqslant 1}}\left(\log f(p^\alpha)-\log f(p^{\alpha-1})\right)\left[\frac{n}{p^\alpha}\right].
\]
Consequently,
\begin{equation}\label{nGn-key}
n\log G_f(n)=\sum_{p\leqslant n}\left[\frac{n}{p}\right]\log f(p)+\sum_{\substack{p^\alpha\leqslant n\\ \alpha\geqslant 2}}\left[\frac{n}{p^\alpha}\right]\log\frac{f(p^\alpha)}{f(p^{\alpha-1})}.
\end{equation}
In particular, if we further assume that $f$ is strongly multiplicative, then 
\begin{equation}\label{nGn-key-strong}
n\log G_f(n)=\sum_{p\leqslant n}\left[\frac{n}{p}\right]\log f(p).
\end{equation}
With regards to this case, we prove the following general result.

\begin{theorem}\label{strong-mult-gem-mean-thm} Let $f$ be a positive strongly multiplicative function such that for any prime $p$ it satisfies $f(p)=\alpha(d)\,p^d+O(p^{d-\delta})$ with $d$ real and $\alpha(d),\delta>0$. Then for any positive integer $r$, there exist computable constants $c_1,\dots,c_r$ such that
\[
G_f(n)=\alpha(d)^M\mathrm{e}^{d(\gamma+E-1)}\varrho_f\,n^d\left(\log n\right)^{\log\alpha(d)}\left(1+\sum_{j=1}^r \frac{c_{j}}{\log^{j}n}+O\Big(\frac{1}{\log^{r+1}n}\Big)\right),
\]
where 
\[
M=\lim_{x\to\infty}\sum_{p\leqslant x}\frac{1}{p}-\log\log x
\] 
is the Meissel--Mertens constant, $\gamma$ is Euler's constant, $E$ is the constant in Mertens' approximation defined by 
\[
E=\lim_{x\to\infty}\sum_{p\leqslant x}\frac{\log p}{p}-\log x,
\]
and $\varrho_f$ is a constant depending on $f$, given by the following product running over primes
\[
\varrho_f=\prod_p\left(\frac{f(p)}{\alpha(d)\,p^d}\right)^{\frac{1}{p}}.
\]
\end{theorem}

As an example satisfying the conditions of Theorem \ref{strong-mult-gem-mean-thm}, we consider the square-free kernel of $n$ defined by $\kappa(n)=\prod_{p|n}p$. 

\begin{corollary}\label{sq-ker-gem-mean-cor}
For any positive integer $r$, there exist computable constants $c_1, \dotsc, c_r$ such that
\[
G_\kappa(n)=\mathrm{e}^{\gamma+E-1}n+\sum_{j=1}^r c_{j}\frac{n}{\log^{j}n}+O\Big(\frac{n}{\log^{r+1}n}\Big).
\]
\end{corollary}

The proof of Theorem \ref{strong-mult-gem-mean-thm} depends on the approximation of the sum on the right hand side of \eqref{nGn-key-strong}. The following key result gives the required approximation, and implies Theorem \ref{strong-mult-gem-mean-thm} immediately.

\begin{theorem}\label{sum-Q-thm} Let $Q(x)=\alpha(d)\,x^d+\mathcal{E}(x)$ with $d$ real and $\alpha(d)>0$, and $\mathcal{E}(x)=O(x^{d-\delta})$ for some fixed $\delta>0$. Moreover, we assume that $Q(n)>0$ for any positive integer $n$. Given any positive integer $r$, there exist computable constants $\eta_0,\eta_1,\dots,\eta_r$ such that
\begin{multline*}
\sum_{p\leqslant n}\left[\frac{n}{p}\right]\log Q(p)
=d\,n\log n+\left(\log\alpha(d)\right)n\log\log n
\\+\eta_0n+\sum_{j=1}^r\eta_j\frac{n}{\log^j n}+O\Big(\frac{n}{\log^{r+1}n}\Big).
\end{multline*}
More precisely 
\begin{equation}\label{eta-0}
\eta_0=M\log\alpha(d)+d(\gamma+E-1)+C_Q,
\end{equation}
where $C_Q$ is an absolute constant in terms of $Q$ defined by
\begin{equation}\label{CQ}
C_Q=\sum_p\frac{1}{p}\log\frac{Q(p)}{\alpha(d)\,p^d}.
\end{equation}
\end{theorem}
\section{Proofs}\label{sec-proofs-props}

We have divided the proof of Theorem \ref{sum-Q-thm} into a sequence of propositions. Let us break up the sum addressed in this theorem as follows
\[
\sum_{p\leqslant n}\left[\frac{n}{p}\right]\log Q(p)
=(\log\alpha(d))\,S_1(n)
+d\,S_2(n)
+S_3(n),
\]
where
\[
S_1(n)=\sum_{p\leqslant n}\left[\frac{n}{p}\right],\quad
S_2(n)=\sum_{p\leqslant n}\left[\frac{n}{p}\right]\log p,\quad
S_3(n)=\sum_{p\leqslant n}\left[\frac{n}{p}\right]\log\left(1+\frac{\mathcal{E}(p)}{\alpha(d) p^d}\right).
\]
An approximation of $S_1(n)$ is well-known due to the notion of the omega function 
\begin{equation}\label{omega}
\omega(k)=\sum_{p|k}1, 
\end{equation}
which counts the number of distinct prime divisors of the positive integer $k$. We observe that
\[
\sum_{k\leqslant n}\omega(k)
=\sum_{k\leqslant n}\sum_{p|k}1
=\sum_{p\leqslant n}\sum_{\substack{k\leqslant n\\ p|k}}1
=\sum_{p\leqslant n}\left[\frac{n}{p}\right].
\]
In 1970 Saffari \cite{saffari} used Dirichlet's hyperbola method to prove 
\begin{equation}\label{Ro-approx-Diac}
\frac{1}{n}\sum_{k\leqslant n}\omega(k)=\log\log n+M+\sum_{j=1}^r\frac{a_j}{\log^j n}+O\Big(\frac{1}{\log^{r+1}n}\Big),
\end{equation}
for each integer $r\geqslant 1$, with 
\begin{equation}\label{aj}
a_j=-\int_1^\infty\frac{\{t\}}{t^2}(\log t)^{j-1}\mathrm{d} t.
\end{equation}
More precisely, it is known \cite{finch} that $a_1=\gamma-1$. Hence, for any positive integer $r$ we obtain
\begin{equation}\label{S1(n)-expan}
S_1(n)=n\log\log n+Mn+\sum_{j=1}^r a_j\frac{n}{\log^j n}+O\Big(\frac{n}{\log^{r+1}n}\Big).
\end{equation}
We mention that later in 1976 Diaconis \cite{diaconis} reproved \eqref{Ro-approx-Diac} by applying Perron's formula on the Dirichlet series $\sum_{n=1}^\infty\omega(n)n^{-s}$ and using complex integration methods. Approximations of $S_2(n)$ and $S_3(n)$ are given in the following propositions.

\begin{proposition}\label{S2(n)-expan-prop} Given any positive integer $r$, there exist computable constants $c_1,\dots,c_r$ such that
\begin{equation}\label{S2(n)-expan}
S_2(n)=n\log n+\left(\gamma+E-1\right)n+\sum_{j=1}^r c_{j}\frac{n}{\log^{j}n}+O\Big(\frac{n}{\log^{r+1}n}\Big).
\end{equation}
\end{proposition}

\begin{proposition}\label{S3(n)-expan-prop} For any fixed $\delta>0$,
\begin{equation}\label{S3(n)-expan}
S_3(n)=C_Qn+O\left(n^{1-\delta}+1+[\delta^{-1}]\log\log n\right).
\end{equation}
\end{proposition}

Note that the big-O term in \eqref{S3(n)-expan} is finally $O(\frac{n}{\log^{r+1}n})$ for each $r>0$. Thus, combining the asymptotic expansions \eqref{S1(n)-expan}, \eqref{S2(n)-expan} and \eqref{S3(n)-expan} completes the proof of Theorem \ref{sum-Q-thm}.

\begin{proof}[Proof of Proposition \ref{S2(n)-expan-prop}]
For each integer $n\geqslant 2$, let 
\[
\lambda(n)=\frac{\log n}{\log\kappa(n)}
\] 
be the index of composition of $n$. In 2005 De Koninck and K\'{a}tai \cite[Theorem 3]{dek-katai} proved that given any positive integer $r$, there exist computable constants $d_1,\dots,d_r$ such that
\begin{equation}\label{sum-1-by-lambda}
U(x):=\sum_{k\leqslant x}\frac{1}{\lambda(k)}=x+\sum_{j=1}^r d_j\frac{x}{\log^j x}+O\Big(\frac{x}{\log^{r+1}x}\Big).
\end{equation}
The above asymptotic expansion is very useful to obtain an asymptotic expansion for $S_2(n)$. Indeed, we observe that
\[
\sum_{k=1}^n\log\kappa(k)
=\sum_{k=1}^n\log\prod_{p|k}p
=\sum_{k=1}^n\sum_{p|k}\log p
=\sum_{p\leqslant n}\left[\frac{n}{p}\right]\log p=S_2(n).
\]
Hence, by Abel summation we get
\[
S_2(n)
=\sum_{k=1}^n\log\kappa(k)
=\sum_{k=2}^n\frac{1}{\lambda(k)}\log k
=U(n)\log n-U(2^-)\log 2-\int_2^n\frac{U(t)}{t}\mathrm{d} t.
\]
To deal with the last integral, we study the functions $\mathrm{L}_j(t)$ defined for each integer $j\geqslant 1$ by the following anti-derivative
\[
\mathrm{L}_j(t):=\int\frac{\mathrm{d} t}{\log^j t}.
\]
Note that $\mathrm{L}_1(t)$ is the logarithmic integral function, which admits the following expansion
\begin{equation}\label{li-expansion}
\mathrm{L}_1(t)=\mathrm{li}(t)=\sum_{i=1}^{r}(i-1)!\,\frac{t}{\log^i t}+O\Big(\frac{t}{\log^{r+1}t}\Big).
\end{equation}
Integrating by parts gives
\[
\mathrm{L}_{j-1}(t)=\int\Big(\frac{1}{\log^{j-1}t}\Big)\left(\mathrm{d} t\right)=\frac{t}{\log^{j-1}t}+(j-1)\int\frac{\mathrm{d} t}{\log^j t}.
\]
Hence, for $j\geqslant 2$ the functions $\mathrm{L}_j(t)$ satisfy the recurrence 
\[
\mathrm{L}_j(t)=\frac{1}{j-1}\mathrm{L}_{j-1}(t)-\frac{t}{(j-1)\log^{j-1}t}.
\]
By repeatedly using this recurrence we deduce that
\[
(j-1)!\,\mathrm{L}_j(t)=\mathrm{li}(t)-\sum_{i=1}^{j-1}(i-1)!\,\frac{t}{\log^i t}.
\]
Hence, by using the expansion \eqref{li-expansion}, for $1\leqslant j\leqslant r$ we obtain
\begin{equation}\label{Ljt-expansion}
\mathrm{L}_j(t)=\sum_{i=j}^{r}\frac{(i-1)!}{(j-1)!}\frac{t}{\log^i t}+O\Big(\frac{t}{\log^{r+1}t}\Big).
\end{equation}
We deduce from the expansion \eqref{sum-1-by-lambda} that
\begin{align*}
\int_2^n\frac{U(t)}{t}\mathrm{d} t
&=\int_2^n\left(1+\sum_{j=1}^r d_j\frac{1}{\log^j t}+O\Big(\frac{1}{\log^{r+1}t}\Big)\right)\mathrm{d} t\\
&=n+\sum_{j=1}^r d_j\mathrm{L}_j(n)-\left(2+\sum_{j=1}^r d_j\mathrm{L}_j(2)\right)+O\Big(\frac{n}{\log^{r+1}n}\Big).
\end{align*}
With $r$ replaced by $r+1$ in \eqref{sum-1-by-lambda}, we obtain
\[
U(n)\log n=n\log n+d_1n+\sum_{j=1}^{r}d_{j+1}\frac{n}{\log^{j}n}+O\Big(\frac{n}{\log^{r+1}n}\Big).
\]
Combining the above expansions yields that
\[
S_2(n)=n\log n+(d_1-1)n+\sum_{j=1}^r\left(d_{j+1}\frac{n}{\log^{j}n}-d_j\mathrm{L}_j(n)\right)-C_r+O\Big(\frac{n}{\log^{r+1}n}\Big),
\]
where $C_r=2+ U(2^-)\log 2+\sum_{j=1}^r d_j\mathrm{L}_j(2)=O_r(1)$. Note that
\begin{multline*}
\sum_{j=1}^r\left(d_{j+1}\frac{n}{\log^{j}n}-d_j\mathrm{L}_j(n)\right)\\
=\sum_{j=1}^r\left(d_{j+1}\frac{n}{\log^{j}n}-\sum_{i=j}^{r}d_j\frac{(i-1)!}{(j-1)!}\frac{n}{\log^i n}\right)+O\Big(\frac{n}{\log^{r+1}n}\Big).
\end{multline*}
Also, an easy computation shows that
\[
\sum_{j=1}^r\left(d_{j+1}\frac{n}{\log^{j}n}-\sum_{i=j}^{r}d_j\frac{(i-1)!}{(j-1)!}\frac{n}{\log^i n}\right)=\sum_{j=1}^r c_{j}\frac{n}{\log^{j}n}+O\Big(\frac{n}{\log^{r+1}n}\Big),
\]
for some computable constants $c_j$ in terms of the $d_j$s. Furthermore we let $c_0=d_1-1$. Consequently,
\[
S_2(n)=n\log n+c_0n+\sum_{j=1}^r c_{j}\frac{n}{\log^{j}n}+O\Big(\frac{n}{\log^{r+1}n}\Big).
\]
It remains to compute the value of $c_0$. To do this, we let
\begin{equation}\label{SMR}
S_2(n)=n\,M(n)-R(n),
\end{equation}
where
\[
M(x):=\sum_{p\leqslant x}\frac{\log p}{p},\qquad\text{and}\qquad
R(n):=\sum_{p\leqslant n}\left\{\frac{n}{p}\right\}\log p.
\]
Theorem 6 of \cite{rs-ille} asserts validity of the double sided inequality
\begin{equation}\label{sum-log-p-by-p}
\log x+E-\frac{1}{2\log x}< M(x)<\log x+E+\frac{1}{2\log x},
\end{equation}
the left hand side for $x>1$ and the right hand side for $x\geqslant 319$. We estimate $R(n)$. Let 
$$
F_1(x)=\sum_{p\leqslant x}\left\{\frac{x}{p}\right\}
\qquad\text{and}\qquad
F_2(x)=\sum_{p^\alpha\leqslant x}\left\{\frac{x}{p^\alpha}\right\}.
$$
Lemma 3 of \cite{lee} asserts that
\[
F_2(n)=(1-\gamma)\frac{n}{\log n}+O\Big(\frac{n}{\log^2 n}\Big).
\]
Note that
\begin{align*}\label{F2-F1-approx}
\F_2(n)-\F_1(n)=\sum_{\substack{p^\alpha\leqslant n\\ \alpha\geqslant 2}}\left\{\frac{n}{p^\alpha}\right\}
&<\sum_{\substack{p^\alpha\leqslant n\\ \alpha\geqslant 2}}1
=\sum_{\substack{p\leqslant n^{\frac{1}{\alpha}}\\ \alpha\geqslant 2}}1
=\sum_{2\leqslant\alpha\leqslant\frac{\log n}{\log 2}}\pi(n^{\frac{1}{\alpha}})\\
&\ll\sum_{2\leqslant\alpha\leqslant\frac{\log n}{\log 2}}\frac{n^{\frac{1}{\alpha}}}{\log n^{\frac{1}{\alpha}}}\leqslant
\frac{n^{\frac12}}{\log n}\sum_{2\leqslant\alpha\leqslant\frac{\log n}{\log 2}}\alpha\ll \sqrt{n}\log n.\nonumber
\end{align*}
Hence
\[
\F_1(n)=(1-\gamma)\frac{n}{\log n}+O\Big(\frac{n}{\log^2 n}\Big).
\]
Let $\varpi(k)$ to be $1$ when $k$ is prime and $0$ otherwise. Abel summation allows us to write
\begin{align*}
R(n)
&=\sum_{k=2}^n\left\{\frac{n}{k}\right\}\varpi(k)\log k\\
&=F_1(n)\log n-F_1(2^-)\log 2-\int_2^n\left(\sum_{p\leqslant t}\left\{\frac{n}{p}\right\}\right)\frac{\mathrm{d} t}{t}\\
&=(1-\gamma)n+O\Big(\frac{n}{\log n}\Big)-\int_2^nO\Big(\frac{t}{\log t}\Big)\frac{\mathrm{d} t}{t}=(1-\gamma)n+O\Big(\frac{n}{\log n}\Big).
\end{align*}
By using \eqref{sum-log-p-by-p} and \eqref{SMR}, we obtain
\[
S_2(n)=n\log n+\left(\gamma+E-1\right)n+O\Big(\frac{n}{\log n}\Big),
\]
implying that $c_0=\gamma+E-1$ and also $d_1=\gamma+E$. 
\end{proof}

\begin{proof}[Proof of Proposition \ref{S3(n)-expan-prop}] We have
\begin{align*}
S_3(n)
&=\sum_{p\leqslant n}\left[\frac{n}{p}\right]\log\frac{Q(p)}{\alpha(d)\,p^d}\\
&=n\sum_p\frac{1}{p}\log\frac{Q(p)}{\alpha(d)\,p^d}-n\sum_{p>n}\frac{1}{p}\log\frac{Q(p)}{\alpha(d)\,p^d}-\sum_{p\leqslant n}\left\{\frac{n}{p}\right\}\log\frac{Q(p)}{\alpha(d)\,p^d}.
\end{align*}
Note that 
\[
\log\frac{Q(p)}{\alpha(d)\,p^d}
=\log\left(1+\frac{\mathcal{E}(p)}{\alpha(d)\,p^d}\right)
=\log\left(1+O\Big(\frac{1}{p^\delta}\Big)\right)=O\Big(\frac{1}{p^\delta}\Big).
\]
Hence
\[
\sum_p\frac{1}{p}\log\frac{Q(p)}{\alpha(d)\,p^d}\ll\sum_p\frac{1}{p^{1+\delta}},
\]
and this implies that in \eqref{CQ} the series defining $C_Q$ is absolutely convergent. Moreover
\[
n\sum_{p>n}\frac{1}{p}\log\frac{Q(p)}{\alpha(d)\,p^d}\ll n\sum_{p>n}\frac{1}{p^{1+\delta}}\ll n\int_n^\infty\frac{\mathrm{d} t}{t^{1+\delta}}\ll n^{1-\delta}.
\]
Also, note that
\[
\sum_{p\leqslant n}\left\{\frac{n}{p}\right\}\log\frac{Q(p)}{\alpha(d)\,p^d}\ll\sum_{p\leqslant n}\frac{1}{p^\delta}.
\]
Hence, if $\delta=1$, then
\[
\sum_{p\leqslant n}\left\{\frac{n}{p}\right\}\log\frac{Q(p)}{\alpha(d)\,p^d}\ll\log\log n,
\]
and if $\delta\neq 1$, then
\[
\sum_{p\leqslant n}\left\{\frac{n}{p}\right\}\log\frac{Q(p)}{\alpha(d)\,p^d}\ll\sum_{p\leqslant n}\frac{1}{p^\delta}
\ll\int_2^n\frac{\mathrm{d} t}{t^{\delta}}\ll n^{1-\delta}.
\]
Combining the above approximations we get \eqref{S3(n)-expan}.
\end{proof}

\begin{proof}[Proof of Theorem \ref{strong-mult-gem-mean-thm}] We observe that for any positive integer $r$
\begin{equation}\label{exp-sum}
\exp\left(\sum_{j=1}^r \frac{c_{j}}{\log^{j}n}+O\Big(\frac{1}{\log^{r+1}n}\Big)\right)=1+\sum_{j=1}^r \frac{c_{j}}{\log^{j}n}+O\Big(\frac{1}{\log^{r+1}n}\Big),
\end{equation}
where the computable constants $c_1,\dots,c_r$ are not necessarily the same on both sides. The relation \eqref{nGn-key-strong} and Theorem \ref{sum-Q-thm} give
\[
\log G_f(n)=d\log n+\left(\log\alpha(d)\right)\log\log n
+\eta_0+\sum_{j=1}^r\frac{\eta_j}{\log^j n}+O\Big(\frac{1}{\log^{r+1}n}\Big).
\]
Taking exponents and using \eqref{exp-sum} completes the proof.
\end{proof}

\begin{proof}[Proof of Corollary \ref{sq-ker-gem-mean-cor}] By using \eqref{nGn-key-strong} we observe that $n\log G_\kappa(n)=S_2(n)$. Thus
\[
\log G_\kappa(n)=\frac{S_2(n)}{n}=\log n+\left(\gamma+E-1\right)+\sum_{j=1}^r \frac{c_{j}}{\log^{j}n}+O\Big(\frac{1}{\log^{r+1}n}\Big).
\]
Taking exponents and considering \eqref{exp-sum} completes the proof.
\end{proof}


\EditInfo{August 09, 2020}{October 05, 2020}{Karl Dilcher}

\end{paper}